# A Conditional Singular Value Decomposition

Qi Deng[1,2,3,4*]

## Abstract

We propose a Conditional Singular Value Decomposition in the form of $A_{\{mn\}} = H_{\{mk\}} B_{\{kl\}} M^*_{\{ln\}}$ for given general matrices $A_{\{mn\}}$ and $B_{\{kl\}}$.

Funding Source: The work was supported by Hubei University of Automotive Technology [grant number BK202209] and Hubei Provincial Bureau of Science and Technology [grant number 2023EHA018].

1. College of Artificial Intelligence, Hubei University of Automotive Technology, Shiyan, China
2. Jack Welch College of Business & Technology, Sacred Heart University, Fairfield, CT, USA
3. School of Accounting, Economics and Finance, University of Portsmouth, Portsmouth, UK
4. Cofintelligence Financial Technology Ltd., Hong Kong and Shanghai, China
*. Corresponding author: dq@huat.edu.cn; dengq@sacredheart.edu

## 1. Introduction

The Singular Value Decomposition (SVD) given as $A_{\{mn\}} = U_{A_{\{mm\}}} \Sigma_{A_{\{mn\}}} V^*_{A_{\{nn\}}}$ is unconditional. However, in financial econometrics applications, there is need to decompose the covariance matrix of one multivariate time series with that of another. To generalize, for two given matrices of $A_{\{mn\}}$ and $B_{\{kl\}}$, under certain conditions, there should exist at least one conditional decomposition that satisfies $A_{\{mn\}} = H_{\{mk\}} B_{\{kl\}} M^*_{\{ln\}}$. In searching for matrix decomposition and factorization literature [e.g., 1-5], we find no direct methodology that addresses this seemingly trivial problem. We propose a conditional decomposition based on the SVD, which we name the "Conditional Singular Value Decomposition" or "Conditional SVR" as a convenient designation.

## 2. Conditional Singular Value Decomposition

For $A_{\{mn\}} \in \mathbb{C}^{m \times n}$ and $B_{\{kl\}} \in \mathbb{C}^{k \times l}$, both $A_{\{mn\}}$ and $B_{\{kl\}}$ have the SVD decompositions as:

$$A_{\{mn\}} = U_{A_{\{mm\}}} \Sigma_{A_{\{mn\}}} V^*_{A_{\{nn\}}} \tag{1a}$$

$$B_{\{kl\}} = U_{B_{\{kk\}}} \Sigma_{B_{\{kl\}}} V^*_{B_{\{ll\}}} \tag{1b}$$

*where*:
1) *$U's$ and $V's$ are square complex unitary matrices;*
2) *$\Sigma's$ are rectangular diagonal matrices with non-negative real numbers on the diagonal.*

**Lemma 1**: There exists a decomposition between $\Sigma_{A_{\{mn\}}}$ and $\Sigma_{B_{\{kl\}}}$ in Equations 1a and 1b, respectively, that satisfies the following:

$$\Sigma_{A_{\{mn\}}} = R_{\{mk\}} \Sigma_{B_{\{kl\}}} S^*_{\{ln\}} \tag{2}$$

*where: $R_{\{mk\}}$ and $S^*_{\{ln\}}$ are diagonal matrices*

1) When $k \geq l$, if $m, n \geq l$, and if the diagonal matrix of $\Sigma_{B_{\{kl\}}}$ has non-zero real numbers on the diagonal, we have a definitely defined diagonal $D_{\{ll\}}$ in:

$$R_{\{mk\}} = \begin{bmatrix} D_{\{ll\}} & 0_{\{l(k-l)\}} \\ 0_{\{(m-l)l\}} & 0_{\{(m-l)(k-l)\}} \end{bmatrix}_{\{mk\}} \tag{3a}$$

$$S^*_{\{ln\}} = [D_{\{ll\}} \quad 0_{\{l(n-l)\}}]_{\{ln\}} \tag{3b}$$

2) When $k \leq l$, if $m, n \geq k$, and if the diagonal matrix of $\Sigma_{B_{\{kl\}}}$ has non-zero real numbers on the diagonal, we have a definitely defined diagonal $D_{\{kk\}}$ in:

$$R_{\{mk\}} = \begin{bmatrix} D_{\{kk\}} \\ 0_{\{(m-k)k\}} \end{bmatrix}_{\{mk\}} \tag{4a}$$

$$S^*_{\{ln\}} = \begin{bmatrix} D_{\{kk\}} & 0_{\{k(n-k)\}} \\ 0_{\{(l-k)k\}} & 0_{\{(l-k)(n-k)\}} \end{bmatrix}_{\{ln\}} \tag{4b}$$

**Proof**:

1) When $k \geq l$ and $m, n \geq l$ holds, we have:

$$\Sigma_{A_{\{mn\}}} = R_{\{mk\}} \Sigma_{B_{\{kl\}}} S^*_{\{ln\}}$$

$$= \begin{bmatrix} D_{\{ll\}} & 0_{\{l(k-l)\}} \\ 0_{\{(m-l)l\}} & 0_{\{(m-l)(k-l)\}} \end{bmatrix}_{\{mk\}} \begin{bmatrix} \Sigma_{B_{\{ll\}}} \\ 0_{\{(k-l)l\}} \end{bmatrix}_{\{kl\}} [D_{\{ll\}} \quad 0_{\{l(n-l)\}}]_{\{ln\}}$$

$$= \begin{bmatrix} D_{\{ll\}} \Sigma_{B_{\{ll\}}} \\ 0_{\{(m-l)l\}} \end{bmatrix}_{\{ml\}} [D_{\{ll\}} \quad 0_{\{l(n-l)\}}]_{\{ln\}} = \begin{bmatrix} D_{\{ll\}} \Sigma_{B_{\{ll\}}} D_{\{ll\}} & 0_{\{l(n-l)\}} \\ 0_{\{(m-l)l\}} & 0_{\{(m-l)(n-l)\}} \end{bmatrix}_{\{mn\}}$$

$$\Rightarrow \Sigma_{A_{\{mn\}}} = \begin{bmatrix} D_{\{ll\}} \Sigma_{B_{\{ll\}}} D_{\{ll\}} & 0_{\{l(n-l)\}} \\ 0_{\{(m-l)l\}} & 0_{\{(m-l)(n-l)\}} \end{bmatrix}_{\{mn\}} \tag{5a}$$

The top-left $l \times l$ sub-diagonal matrix of $\Sigma_{A_{\{mn\}}}$, or $\Sigma_{A_{\{ll\}}}$, exists and therefore:

$$\Sigma_{A_{\{ll\}}} = D_{\{ll\}} \Sigma_{B_{\{ll\}}} D_{\{ll\}} = D_{\{ll\}} D_{\{ll\}} \Sigma_{B_{\{ll\}}} \tag{5b}$$

$$D_{\{ll\}} = \left( \Sigma_{A_{\{ll\}}} \Sigma_{B_{\{ll\}}}^{-1} \right)^{\frac{1}{2}} \tag{5c}$$

Equation 5c holds only if $\Sigma_{B_{\{ll\}}}$ is inversible ($\Sigma_{B_{\{ll\}}}^{-1}$ exists), therefore we prove Lemma 1 when $k \geq l$ and $m, n \geq l$, by tightening the condition that the top $l \times l$ sub-diagonal matrix of $\Sigma_{B_{\{kl\}}}$, or $\Sigma_{B_{\{ll\}}}$, must have non-zero diagonal elements.

2) When $k \leq l$ and $m, n \geq k$ holds, we have:

$$\Sigma_{A_{\{mn\}}} = R_{\{mk\}} \Sigma_{B_{\{kl\}}} S^*_{\{ln\}}$$

$$= \begin{bmatrix} D_{\{kk\}} \\ 0_{\{(m-k)k\}} \end{bmatrix}_{\{mk\}} \begin{bmatrix} \Sigma_{B_{\{kk\}}} & 0_{\{k(l-k)\}} \end{bmatrix}_{\{kl\}} \begin{bmatrix} D_{\{kk\}} & 0_{\{k(n-k)\}} \\ 0_{\{(l-k)k\}} & 0_{\{(l-k)(n-k)\}} \end{bmatrix}_{\{ln\}}$$

$$= \begin{bmatrix} D_{\{kk\}}\Sigma_{B_{\{kk\}}} & 0_{\{k(l-k)\}} \\ 0_{\{(m-k)k\}} & 0_{\{(m-k)(l-k)\}} \end{bmatrix}_{\{ml\}} \begin{bmatrix} D_{\{kk\}} & 0_{\{k(n-k)\}} \\ 0_{\{(l-k)k\}} & 0_{\{(l-k)(n-k)\}} \end{bmatrix}_{\{ln\}}$$

$$= \begin{bmatrix} D_{\{kk\}}\Sigma_{B_{\{kk\}}}D_{\{kk\}} & 0_{\{k(n-k)\}} \\ 0_{\{(m-k)k\}} & 0_{\{(m-k)(n-k)\}} \end{bmatrix}_{\{mn\}}$$

$$\Rightarrow \Sigma_{A_{\{mn\}}} = \begin{bmatrix} D_{\{kk\}}\Sigma_{B_{\{kk\}}}D_{\{kk\}} & 0_{\{k(n-k)\}} \\ 0_{\{(m-k)k\}} & 0_{\{(m-k)(n-k)\}} \end{bmatrix}_{\{mn\}} \tag{6a}$$

The top-left $k \times k$ sub-diagonal matrix of $\Sigma_{A_{\{mn\}}}$, or $\Sigma_{A_{\{kk\}}}$ exists and therefore:

$$\Sigma_{A_{\{kk\}}} = D_{\{kk\}}\Sigma_{B_{\{kk\}}}D_{\{kk\}} = D_{\{kk\}}D_{\{kk\}}\Sigma_{B_{\{kk\}}} \tag{6b}$$

$$D_{\{kk\}} = \left(\Sigma_{A_{\{kk\}}}\Sigma_{B_{\{kk\}}}^{-1}\right)^{\frac{1}{2}} \tag{6c}$$

Equation 6c holds only if $\Sigma_{B_{\{kk\}}}$ is inversible ($\Sigma_{B_{\{kk\}}}^{-1}$ exists), therefore we prove Lemma 1 when $k \leq l$ and $m, n \geq k$, by tightening the condition that the left $k \times k$ sub-diagonal matrix of $\Sigma_{B_{\{kl\}}}$, or $\Sigma_{B_{\{kk\}}}$, must have non-zero diagonal elements.

**Theorem 1**: Let $A_{\{mn\}} \in \mathbb{C}^{m \times n}$ and $B_{\{kl\}} \in \mathbb{C}^{k \times l}$, there exists at least one conditional decomposition that satisfies the following:

$$A_{\{mn\}} = H_{\{mk\}}B_{\{kl\}}M_{\{ln\}}^* \tag{7}$$

Under either of the following conditions:

Condition 1) When $k \geq l$, if $m, n \geq l$, and if the diagonal matrix of $\Sigma_{B_{\{kl\}}}$ in Equation 1b has non-zero real numbers on the diagonal;

Condition 2) When $k \leq l$, if $m, n \geq k$, and if the diagonal matrix of $\Sigma_{B_{\{kl\}}}$ in Equation 1b has non-zero real numbers on the diagonal.

**Proof**:

We first substitute the $B_{\{kl\}}$ term in Equation 7 with the RHS of Equation 1b, resulting:

$$A_{\{mn\}} = H_{\{mk\}} \left( U_{B_{\{kk\}}} \Sigma_{B_{\{kl\}}} V^*_{B_{\{ll\}}} \right) M^*_{\{ln\}} = \left( H_{\{mk\}} U_{B_{\{kk\}}} \right)_{\{mk\}} \Sigma_{B_{\{kl\}}} \left( M_{\{nl\}} V_{B_{\{ll\}}} \right)^*_{\{ln\}} \quad (8a)$$

We then substitute the $\Sigma_{A_{\{mn\}}}$ term in Equation 1a with the RHS of Equation 2, resulting:

$$A_{\{mn\}} = U_{A_{\{mm\}}} \left( R_{\{mk\}} \Sigma_{B_{\{kl\}}} S^*_{\{ln\}} \right) V^*_{A_{\{nn\}}} = \left( U_{A_{\{mm\}}} R_{\{mk\}} \right)_{\{mk\}} \Sigma_{B_{\{kl\}}} \left( V_{A_{\{nn\}}} S_{\{nl\}} \right)^*_{\{ln\}} \quad (8b)$$

Comparing the RHSs of Equations 8a and 8b, we get the follows:

$$\left( H_{\{mk\}} U_{B_{\{kk\}}} \right)_{\{mk\}} = \left( U_{A_{\{mm\}}} R_{\{mk\}} \right)_{\{mk\}} \quad (9a)$$

$$\left( M_{\{nl\}} V_{B_{\{ll\}}} \right)^*_{\{ln\}} = \left( V_{A_{\{nn\}}} S_{\{nl\}} \right)^*_{\{ln\}} \Rightarrow \left( M_{\{nl\}} V_{B_{\{ll\}}} \right)_{\{nl\}} = \left( V_{A_{\{nn\}}} S_{\{nl\}} \right)_{\{nl\}} \quad (9b)$$

Solving Equations 9a and 9b, we get:

$$H_{\{mk\}} = U_{A_{\{mm\}}} R_{\{mk\}} U^*_{B_{\{kk\}}} \quad (10a)$$

$$M_{\{nl\}} = V_{A_{\{nn\}}} S_{\{nl\}} V^*_{B_{\{ll\}}} \quad (10b)$$

If Condition 1 holds, we substitute $R_{\{mk\}}$ and $S_{\{nl\}}$ in Equations 10a and 10b with Equations 3a and 3b, respectively, with $D_{\{ll\}}$ given by Equation 5C. If Condition 2 holds, we substitute $R_{\{mk\}}$ and $S_{\{nl\}}$ with Equations 4a and 4b, respectively, with $D_{\{kk\}}$ given by Equation 6C. As all the matrices exist, we thus prove Theorem 1.

It is trivial to prove that both $HH^*$ and $MM^*$ are symmetric matrices for any permutation of $m, n, k, l$ that satisfies the constraints in either condition. For example, from Equation 10a we get:

$$\left( H_{\{mk\}} H^*_{\{km\}} \right)^* = H_{\{mk\}} H^*_{\{km\}}$$

3. **A Special Case**

In case that $m = n = k = l$, Equation 7 is reduced to:

$$A = HBH^* \quad (11)$$

$A$ and $B$ have the SVD decomposition given as:

$$A = U_A \Sigma_A U_A^* \tag{12a}$$

$$B = U_B \Sigma_B U_B^* \tag{12b}$$

    1) U's are square complex unitary matrices;
    2) Σ's are rectangular diagonal matrices with non-negative real numbers on the diagonal.

And there exists a decomposition between $\Sigma_A$ and $\Sigma_B$ as:

$$\Sigma_A = R \Sigma_B R^* \tag{13}$$

    where: R is a diagonal matrix with real numbers on the diagonal

$$\Rightarrow \Sigma_A = RR^* \Sigma_B = RR \Sigma_B$$

$$\Rightarrow R = (\Sigma_A \Sigma_B^{-1})^{\frac{1}{2}} \tag{14}$$

By substituting $\Sigma_A$ in Equation 12a with the RHS of Equation 13 we get:

$$A = U_A (R \Sigma_B R^*) U_A^*$$

$$\Rightarrow A = (U_A R) \Sigma_B (U_A R)^* \tag{15a}$$

Also, substitute $B$ in Equation 11 with the RHS of Equation 12b:

$$A = H(U_B \Sigma_B U_B^*) H^*$$

$$\Rightarrow A = (H U_B) \Sigma_B (H U_B)^* \tag{15b}$$

By comparing Equation 15a and Equation 15b we get:

$$U_A R = H U_B$$

$$\Rightarrow H = U_A R U_B^* \tag{16}$$

Equation 16 solves $H$ in proposition $A = HBH^*$, with $HH^*$ being a symmetric matrix.

## 4. Discussion

In this short note we prove that there exists a Conditional SVD of $A_{\{mn\}} = H_{\{mk\}} B_{\{kl\}} M_{\{ln\}}^*$, and provide an analytical solution for it. We contribute to the literature of matrix decomposition, and especially that of the SVD, by proposing a technique that directly addresses the decomposition of $A_{\{mn\}}$ under the condition that $B_{\{kl\}}$ is also given. We also provide a special case, that when $m = n = k = l$, a reduced conditional SVD in the form of $A = HBH^*$ exists.

**Declaration of competing interest**

There is no competing interest.

**Data availability**

No data was used for the research described in the article.